\documentstyle[epsfig,11pt]{article}
\onecolumn 
\begin{document}


\newcommand {\vare}{\varepsilon}
\newcommand {\grad}{\bigtriangledown}
\newcommand {\RR}{\,\hbox{\rm R}\!\!\!\!\!\hbox{\rm I}\,\,\,}
\newcommand{\square}{\hbox{${\vcenter{\hrule height.4pt
   \hbox{\vrule width.4pt height6pt \kern6pt
     \vrule width.4pt} \hrule height.4pt}}$}}
\newcommand {\ProofEnd} {\hfill \nobreak \square \medbreak}


\newtheorem {theorem} {Theorem}
\newtheorem {proposition} [theorem] {Proposition}
\newtheorem {lemma} [theorem] {Lemma}
\newtheorem {definition} [theorem] {Definition}
\newtheorem {corollary} [theorem] {Corollary}
\newtheorem {remark} [theorem] {Remark}
\newtheorem {note} [theorem] {Note}
\newtheorem {example} [theorem] {Example}

\begin{titlepage} 
\begin{center}
{\huge\bf Harmonic Functions on Manifolds with\\
\vspace*{2ex} Nonnegative Ricci Curvature\\
\vspace*{2ex} and Linear Volume Growth}

\vspace*{7ex} 
Christina Sormani$ ^\dagger$\footnote{\\
        \hspace*{1em}$ ^\dagger$The author can be reached by e-mail at
       {\tt sormani@math.jhu..edu}\\ 
       {\it Department of Mathematics\\
        Johns Hopkins University\\
        3400 North Charles Street\\
        Baltimore, MD 21218 }\\} 
\date{July 1997}
\end{center}
\end{titlepage} 
\pagestyle{plain}  
\pagenumbering{arabic}



\newpage
\vspace{2cm} 


{\bf
In this paper we prove that if a complete noncompact manifold with
nonnegative Ricci curvature and linear volume growth has a nonconstant
harmonic function of at most polynomial growth, then the manifold 
splits isomtrically.}


In 1975, Shing Tung
Yau proved that a complete noncompact manifold with nonnegative
Ricci curvature has no nonconstant harmonic functions of
sublinear growth [Yau2]. That is, if
\begin{equation}
\limsup_{R \to \infty} \frac {\max_{B_p(R)}|f|}{R} =0
\end{equation}
and if $f$ is harmonic, then $f$ is a constant.  In the same
paper, Yau used this result to prove that a complete noncompact
manifold with nonnegative Ricci curvature has at least linear
volume growth,
\begin{equation}
\liminf_{R \to \infty} \frac {Vol(B_p(R))}{R} =C \in (0,\infty].
\end{equation}

There are many manifolds with nonnegative Ricci curvature
that actually have linear volume growth
\begin{equation}
\limsup_{R \to \infty} \frac {Vol(B_p(R))}{R} =V_0 <\infty.
\end{equation}
Some interesting examples of such manifolds can be found in [So1].

In this paper, we prove the following theorem concerning harmonic functions
on these manifolds.  

\begin{theorem}  Let $M^n$ be a complete noncompact manifold with
nonnegative Ricci curvature and 
at most linear volume growth,
\begin{equation}
\limsup_{R\to\infty} \frac {Vol(B_p(R))}{R} =V_0<\infty  
\end{equation}
If there exists a nonconstant harmonic function, $f$, of
polynomial growth of any given degree $q$,
\begin{equation} \label{poly}
 \Delta f = 0 \quad \textrm{  and   }\quad |f(x)| \le C(d(x,p)^q+1),
\end{equation}
then the manifold splits isometrically, $M^n = N^{n-1}\times \RR$.
\end{theorem}

Harmonic functions of polynomial growth have been an object of study for
some time.  Until recently it was not known whether the space of
harmonic functions of polynomial growth of a given degree on a manifold
with nonnegative Ricci curvature was finite dimensional.  
Atsushi Kasue proved this result with various additional assumptions in
[Kas1, Kas2].  Tobias Colding and Bill Minnicozzi have recently proven 
that this space is indeed finite dimensional with no additional assumptions
[CoMin].   With our stronger condition of 
linear volume growth, we are able to prove that this space is 
only one dimensional directly using a gradient estimate of Cheng and Yau 
[ChgYau] and results from [So1, So2].

For background material consult the textbooks [SchYau] and [Li].  

The author would like to thank Professor Shing-Tung Yau for conjecturing this
theorem and for his encouragement during her year at 
Harvard University.  She would also thank to thank Professors Tobias Colding,
Jeff Cheeger, William Minicozzi and Gang Tian for their continued interest 
in her work.

\begin{section}{Background}

\indent
A {\em ray}, $\gamma:[0,\infty)\mapsto M^n$, is a geodesic
which is minimal on any subsegment, $d(\gamma(t),\gamma(s))=|t-s|$.
Every complete noncompact Riemannian manifold contains a ray.  Given
a ray, one can define its associated {\em Busemann function}, 
$b: M^n \to \RR$,
as follows:
\begin{equation}
b(x)=\lim_{R \to \infty} \big( R-d(x,\gamma(R)) \big).
\end{equation}
The Busemann function is a Lipschitz function whose gradient has unit length
almost everywhere.  [Bu].  

In Euclidean space, the level sets of the Busemann function associated
with a given ray  are the planes perpendicular to the
given ray.  In contrast,  the Busemann function
defined on a manifold with nonnegative Ricci curvature and linear volume 
growth has compact level sets with bounded diameter growth [So1, Thm 15].
In that paper, the author also proved that if such a manifold is not an
isometrically split manifold, then the Busemann function is bounded below
and $b^{-1}((-\infty, r])$ is a compact set for all $r$ [So1, Cor 19].

\begin{lemma} 
Let $M^n$ be a complete manifold with nonnegative Ricci 
curvature.  Suppose that there is a Busemann function, $b$, which
is bounded below by $b_{min}$ and that
diameter of the level sets grows at most linearly,
\begin{equation} \label{diam}
diam(b^{-1}(b_{min}+r)) \le C_D(r+1).
\end{equation}
Then there exists a universal constant, $C_n$, depending only on the 
dimension, $n$, such that
any harmonic function, $f$, satisfies the following gradient estimate:
\begin{equation} \label{grad}
\sup_{b^{-1}([b_{min}, b_{min}+r))} \,|\grad f|\,\, 
\le \,\,\frac {C_n}{2(r+D)}\, \sup_{b^{-1}(b_{min}+2(r+D))} \,|f| 
\end{equation} 
for all $D\ge C_D(r+1)$.
\end{lemma}

\noindent {\bf Proof of Lemma:}
First note that the boundary of the compact set, $b^{-1}([b_{min}, r))$,
is just $b^{-1}(r)$.  So,
by the maximum principal, we know that for any harmonic function, $f$,
\begin{equation}  \label{max}
\max_{b^{-1}([b_{min}, r))}\,f \,\le \max_{b^{-1}(r)}\,f 
\,\,\,\textrm{ and}\,\,\,
\min_{b^{-1}([b_{min}, r))}\,f\, \ge \min_{b^{-1}(r)}\,f. 
\end{equation}

Furthermore, Cheng and Yau have proven the following gradient estimate
for harmonic functions on balls in manifolds with nonnegative Ricci curvature,
\begin{equation} \label{yau}
\sup_{B_p(a/2)}\, | \grad f| \,\le\, \frac {C_n}{a}\, \sup_{B_p(a)} \,|f|
\end{equation}
where $C_n$ is a universal constant depending only on the dimension, $n$.
[ChgYau, Thm 6], see also [SchYau, p21, Cor 2.2].  This will be the constant 
in (\ref{grad}).  Thus, we need only relate balls to regions defined by the 
Busemann function to prove the theorem.

Let $x_0$ be a point in $b^{-1}(b_{min})$.  Note that 
\begin{equation} \label{in}
B_{x_0}(a) \subset b^{-1}([b_{min}, b_{min}+a))
\end{equation} 
because the triangle inequality implies that
$$
b(x)=\lim_{R \to \infty} R-d(x, \gamma(R))
    \le \lim_{R \to \infty} R-d(x_0, \gamma(R)) +d(x_0,x) =b(x_0)+d(x_0,x).
$$

On the other hand, using our diameter bound in (\ref{diam}), we claim that 
\begin{equation} \label{out}
b^{-1}([b_{min}, b_{min}+r)) \subset B_{x_0}(r + D) \qquad 
\forall D \ge C_D(r+1).
\end{equation} 
To see this we will construct a ray, $\sigma$, emanating from $x_0$ such that
for all $t \ge 0$,
$\sigma(t) \in b^{-1}(r_{min}+t)$.  
Then, for any $y \in b^{-1}([b_{min},b_{min}+r))$, we let $t=b(y)$ and we have
$$
d(x_0,y) \le d(x_0, \sigma(t)) +d(\sigma(t),y)
\le t + diam(b^{-1}(r_{min}+t))\le r+D. 
$$
which implies (\ref{out}).
The ray, $\sigma$, is constructed by taking a limit of minimal geodesics, 
$\sigma_i$,
from $x_0$ to $\gamma(R_i)$.  A subsequence of such a sequence of minimal 
geodesics always converges.   The limit ray satisfies the required property,
\begin{eqnarray*}
b(\sigma(t))&=&\,\lim_{i\to\infty}\,\, b(\sigma_i(t))
=\lim_{i\to\infty} \lim_{R\to\infty} R-d(\sigma_i(t),\gamma(R))\\
&=&\lim_{R\to\infty} R-(d(\sigma_i(0),\gamma(R))-t)\, =\,\,b(x_0)+t.
\end{eqnarray*}

We can now combine the relationships between Busemann regions
and balls, (\ref{out}) and (\ref{in}),
with the gradient estimate, (\ref{yau}), and the maximum principal,
(\ref{max}),
to prove the lemma.  That is, for all $D \ge C_D(r+1)$, we have
\begin{eqnarray*}
\sup_{b^{-1}([b_{min},b_{min}+r))}\,
|\grad f|& \le& \sup_{B_{x_0}(r+D)}|\grad f| \textrm{  by (\ref{out}), }\\
&\le & \frac {C_n}{2(r+D)} \,\sup_{B_{x_0}(2(r+D))}|f| \\
&\le & \frac {C_n}{2(r+D)} \,\sup_{b^{-1}([b_{min},b_{min}+2(r+D)))}|f|\\
&\le & \frac {C_n}{2(r+D)} \,\sup_{b^{-1}(b_{min}+2(r+D))}|f|.
\end{eqnarray*}
\ProofEnd

We employ this lemma and elements of the proof to prove our theorem.
\end{section}

\begin{section} {Proof of the Theorem}

\indent
The given  manifold, $M^n$, has nonnegative Ricci curvature and linear volume
growth.  We will assume that $M^n$
doesn't split isometrically and demonstrate that the harmonic functions of 
polynomial growth must be constant.  Since the manifold doesn't split
isometrically and has linear volume growth, any Busemann function, $b$, has
a minimum value by [So1, Cor 23].   Furthermore, by [So2, Thm ?],
the diameters of the level sets of the Busemann function grow sublinearly.
Thus we satisfy the hypothesis of Lemma 2 with  $C_D=1$ in (\ref{diam}).  

Let $M(r)= \max_{b^{-1}(b_{min}+r)} |f|$, where $f$ is a harmonic function
of polynomial growth. 
Note that $M$ is an nondecreasing function by the maximum principal,
(\ref{max}).
By the lemma, we know that for all $r \ge b_{min}$ and for all 
$D \ge (r+1)$, we can bound the gradient of $f$ in terms of $M$,
\begin{equation}
\sup_{b^{-1}([b_{min},b_{min} +r))}\,|\grad f|\,\, \le \,\,
\frac {C_n M(2(r+D))}{2(r+D)}.
\end{equation}

Since $b^{-1}(r)$ is compact, there exists $x_r, y_r \in b^{-1}(b_{min}+r)$
such that 
\begin{equation}
f(x_r)=\min_{b^{-1}(b_{min}+r)} f \quad \textrm{  and  }\quad 
f(y_r)=\max_{b^{-1}(b_{min}+r)} f.
\end{equation}  
We claim that, for $r$ sufficiently large, $M(r) \le f(y_r)-f(x_r)$.   

First recall that if $f$ is a positive or negative
harmonic function on a manifold with nonnegative Ricci curvature, then
$f$ must be constant [Yau1, Cor 1].   So there exists a point 
$z \in M^n$ such that $f(z)=0$.  Thus, by the maximum principal,
if $r \ge b(z)$ we know that $f(y_r) \ge 0$ and $f(x_r)\le 0$.  So
$M(r) = max(f(y_r), -f(x_r)) \le f(y_r)-f(x_r)$.

We can now estimate $M(r)$ from above in terms of the gradient of $f$
and the diameter of the level set, $b^{-1}(r)$.  First we join
$x_r$ to $y_r$ by a smooth minimal geodesic, $\gamma_r$.  Note that
the length of $\gamma_r$, is less than or equal to $diam(b^{-1}(r))$
by the definition of diameter.  So
$\gamma_r \subset b^{-1}(b_{min}, r +diam(b^{-1}(r)))$.
Thus for all $r\ge b(z)$, for all $D \ge (r+1)$, we have 
\begin{eqnarray*}
M(r)\le f(y_r)-f(x_r) 
&\le & \int_0^{L(\gamma_r)} \frac d {dt} f(\gamma(t))\,dt \\
&\le&\int_0^{L(\gamma_r)} |\grad f|\, |\gamma'(t)|\,dt \\
&\le& \sup_{b^{-1}([b_{min}, r +diam(b^{-1}(r))))} |\grad f| \,\,\,\,
\int_0^{L(\gamma_r)} |\gamma'(t)|\,dt \\
&\le & \frac {C_n M(2(r+diam(b^{-1}(r)) +D))}{2(r+diam(b^{-1}(r)) +D)}\,\, 
diam(b^{-1}(r))\\
& \le & C_n M(2(r+diam(b^{-1}(r)) +D))\,\,\frac{ diam(b^{-1}(r))}{2r}\\
& \le & C_n M(2(r+(r+1) + D))\,\frac{ diam(b^{-1}(r))}{2r}
\end{eqnarray*}
Setting $D=r+1$ and taking $r \ge 1$, we have
\begin{equation} 
M(r) \le C_n \,M(6r)\,\, \frac{ diam(b^{-1}(r))}{2r}.
\end{equation}

Recall that our manifold has sublinear diameter growth
by [So2, Thm ?].  So, given any $\delta >0$,
we can find $R_\delta \ge 1$ such that
\begin{equation} \label{delta} 
\frac{diam(b^{-1}(r))}{2r}\, <\, \delta  \qquad \forall r\ge R_\delta.
\end{equation}
Then, for all $k \in {\bf N}$ and for all $R \ge R_\delta$, we have
\begin{equation} \label{iterate}
M(R)\,\le\, C_n \,M(6R)\, \delta\, \le 
\textrm{ ... }\le\, C_n^k\, M(6^k R)\, \delta^k.
\end{equation}

Now $f$ has polynomial growth of order $q$, (\ref{poly}), so
\begin{equation}
M(r)\,=\, \max_{x\in b^{-1}(b_{min}+r)}\,|f(x)| 
\,\le\, \max_{x\in b^{-1}\,(b_{min}+r)}C(d(x,x_0)^q+1).
\end{equation}
Applying (\ref{out}) with $C_D=1$ and $D=C_D(r+1)$, we get
\begin{equation}
M(r)\,\le \, C((r+(r+1))^q+1)\,\le\, C (6r^2)^q \qquad \forall \,r \ge 1.
\end{equation}

Substituting this information into (\ref{iterate}), we get
\begin{eqnarray*}
M(R)&\le& C_n^k\, C \,\bigg( 6(6^k R)^2\bigg)^q\, \delta^k \\
& \le & C\, 6^q \, R^{2q}\, (C_n 6^{2q} \delta) ^k 
\qquad \forall R \ge R_\delta.
\end{eqnarray*}.

Fix $\delta = 1/(2C_n 6^{2q})$, so $R_\delta$ is fixed by (\ref{delta}).
Then, for all $R \,\ge \,R_\delta$, we have
\begin{equation}
M(R) \,\le\, \lim_{k \to \infty}\, C\, 6^q\, R^{2q} \,(1/2)^k\, =\,0.
\end{equation} 
Since $M(r)$ is nondecreasing and nonnegative, $M(r)=0$ everywhere.
Thus, $f$ is a constant.

\ProofEnd

\end{section}

\end{document}